\input amstex.tex
\documentstyle{amsppt}
\magnification=\magstep1

\font\sc=cmcsc10

\newcount\t
\newcount\q
\newcount\x
\t=0 \q=0 \x=0

\long\def\se#1{\advance\q by1
\t=0 \x=0
\bigskip
\noindent
\S\number\q \; {\sc #1}\nopagebreak\smallskip\nopagebreak}

\long\def\form#1{\global\advance\t by1%
$$ #1 \tag\number\q.\number\t $$}

\long\def\thm#1{\medskip\noindent
\advance\x by1
{\sc Theorem \number\q.\number\x.}\quad {\sl #1}\smallskip}

\long\def\lem#1{\medskip\noindent
\advance\x by1
{\sc Lemma \number\q.\number\x}\quad {\sl #1}\smallskip}

\def\dimo{\noindent {\sc Proof}\quad}

\long\def\prop#1{\medskip\noindent
\advance \x by1
{\sc Proposition \number\q.\number\x}\quad {\sl #1}\smallskip}

\long\def\cor#1{\medskip\noindent
\advance\x by1
{\sc Corollary \number\q.\number\x}\quad {\sl #1}\smallskip}

%\overfullrule=0pt
%%%%%%%%%%%%%%%%%%%%%%%%%%%%%%%%%%%%%%%%%%%%%%%%%%%%%%%%%%%%%%%%%%%%%%%%%%
\topmatter
\title
On the failure of the Poincar\'e Lemma for 
$\pmb{\bar\partial_M}$\;  II
\endtitle
\author
C.Denson Hill and Mauro Nacinovich
\endauthor
\address C.D.Hill, Department of Mathematics, SUNY at Stony Brook,
Stony Brook N.Y. 11794 - USA \endaddress
\address M.Nacinovich: Dipartimento di Matematica, Universit\`a di Roma
``Tor Vergata'', via della Ricerca Scientifica, 00133 Roma - ITALY
\endaddress
\subjclass 32V05 35N15 \endsubjclass
\keywords Poincar\'e lemma, tangential Cauchy-Riemann complex\endkeywords
\endtopmatter
\document
\se{introduction}
The purpose of this paper is to repair some inaccuracies in the
formulation of the main result of \cite{HN1}.
As were written, the main theorems 4.1 and 4.2 of \cite{HN1}
are in fact in contradiction to earlier results of one of the
authors \cite{N2}. In the process of writing an erratum, we actually
discovered some new phenomena. As we found these quite interesting,
it has lead us to incorporate the needed corrections into this
self contained article. 
\par
An unfortunate misprint, in which $R^{-1}$ got replaced by $R$,
led the authors to misinterpret what their proof in \cite{HN1} actually
demonstrated. Upon closer scrutiny, we realized that there are two
distinct ways to proceed.
\par
     One is that we may still obtain the original conclusions of our
main theorems, provided we slightly strengthen our original hypothesis
(cf. Theorems 5.1 and 5.2). This however entails a much more complicated
argument, involving the $CR$ structure of the characteristic bundle, which
is of considerable independent interest (cf. Theorems 4.3 and 4.6).
\par
     The other is that if we stick to our original hypothesis, then the
conclusions we obtain are slightly weaker than originally claimed, but
in our opinion still interesting. In fact a new invariant comes into
play, which measures the rate of shrinking, even in the situation where
the local Poincar\'e lemma is valid.
\par
      Recall that here we make an important
distinction between the vanishing of the cohomology and the validity of
the Poincar\'e lemma: Consider the inhomogeneous problem
$\bar\partial_M u = f$,
to be solved for $u$, with given data 
$f$ satisfying $\bar\partial_M f = 0$ in some
domain $U$ containing a point $x_0$. 
The vanishing of the cohomology in $U$
refers to the situation in which, no matter how $f$ is prescribed in $U$,
there is always a solution $u$ in $U$ (i.e., no shrinking). The validity
of the Poincar\'e lemma at $x_0$ requires only that a solution $u$ exist in
a smaller domain $V_f$ with $x_0\in V_f\subset U$ (i.e., 
there is some shrinking which might, in principle, depend on $f$). Our
new invariant measures the relative rate of shrinking of $V_f$ with
respect to radius$(U)$, as $U$ shrinks to the point $x_0$.\par
 Under our original
hypothesis, we are able to show that the cohomology of small convex
neighborhoods of $x_0$ is always infinite dimensional, with respect to
any choice of the Riemannian metric (cf. Theorems 7.2 and 7.3). This
means that special shapes are needed, if one is to have the vanishing
of the cohomology for small sets.\par
In \S 6 we have listed a number of natural examples satisfying the 
slightly strengthened hypothesis. They illustrate how common it is for the 
Poincar\'e lemma to fail.
\par
        At the end of \S 7 we give a very simple example
   satisfying our original hypothesis. It illustrates how,
   even if the Poincar\'e lemma were to be valid, a shrinking
   must occur which goes like $\roman{radius}(V_f) \simeq C r^{3/2}$, where
   $\roman{radius}( U)\simeq  r$, as $r@>>> 0$.

     The first author would like to express his appreciation for the
kind hospitality of the Universit\`a di Roma "Tor Vergata" and the
Humboldt Universit\"at zu Berlin, and in particular Professor
J\"urgen Leiterer. The second author would like to express his
thanks to Fabio Nicola from Turin for pointing out to us 
inconsistencies in \cite{HN1}.
%%%%%%%%%%%%%%%%%%%%%%%%%%%%%%%%%%%%%%%%%%%%%%%%%%%%%%%%%%%%%%%%%%%%%%

\se{A priori estimates}\edef\sectfunctanal{\number\q}
In this section we rehearse the facts of functional analysis that
we shall use later for the discussion of the local cohomology of
the $\bar\partial_M$-complex. For the proofs we refer to
\cite{AFN}, \cite{N1} and \cite{HN1}.
Let $M$ be a 
paracompact smooth differentiable manifold, of real
dimension $N$. We consider a sequence
of complex vector bundles $\{E^q@>>>M\}_{q=0,1,\hdots}$
(with $E^q$ of rank $r_q$),
and a complex of linear partial differential operators\,:
\form{\CD
\Cal{C}^{\infty}(M,E^0)@>{A_0}>>\Cal{C}^{\infty}(M,E^1)@>{A_1}>>
\Cal{C}^{\infty}(M,E^2)@>>>
\cdots\endCD}\edef\complesso{\number\q.\number\t}
This means that, for
all $i=0,1,\hdots$\,: 
\roster\item"($i$)" 
$A_i:\Cal{C}^\infty(M,E^q)@>>>\Cal{C}^\infty(M,E^{q+1})$
are linear and $\roman{supp}\left(A_q(u)\right)\subset
\roman{supp}(u)$ for all $u\in\Cal{C}^\infty(M,E^q)$
(Peetre's theorem), 
\item"($ii$)" $A_{q+1}\circ A_q=0$\,
\, for
all $q=0,1,\hdots$.
\endroster
\par
We denote by $\Cal{E}^q$ the sheaf of germs of smooth sections of
$E^q$ and, for every open $U\subset M$, set
$\Cal{E}^q(U)=\Cal{C}^\infty(U,E^q)$.\par
For every open $U\subset M$ and $x_0\in M$, we obtain
complexes
$$\CD
\Cal{E}^0(U)@>{A_0}>>\Cal{E}^1(U)@>{A_1}>>
\Cal{E}^2(U)@>>>\cdots\endCD \tag $U,\Cal{E}^*,A_*$ $$
and
$$\CD
\Cal{E}^0_{x_0}@>{A_0}>>\Cal{E}^1_{x_0}@>{A_1}>>
\Cal{E}^2_{x_0}@>>>\cdots\endCD \tag $(x_0),\Cal{E}^*,A_*$ $$
We denote their cohomology groups by\,:
\form{H^q(U,\Cal{E}^*,A_*)=\dsize\frac{\ker\left(A_q:\Cal{E}^q(U)
@>>>\Cal{E}^{q+1}(U)\right)}{\roman{Image}
\left(A_{q-1}:\Cal{E}^{q-1}(U)
@>>>\Cal{E}^{q}(U)\right)} }
and
\form{
H^q((x_0),\Cal{E}^*,A_*)=\dsize\frac{\ker\left(A_q:\Cal{E}^q_{x_0}
@>>>\Cal{E}^{q+1}_{x_0}\right)}{\roman{Image}
\left(A_{q-1}:\Cal{E}^{q-1}_{x_0}
@>>>\Cal{E}^{q}_{x_0}\right)}
=\varinjlim_{U\ni x_0}H^q(U,\Cal{E}^*,A_*)\,, 
}
respectively. \par
When $H^q((x_0),\Cal{E}^*,A_*)=\{0\}$,
we say that (\complesso) admits the Poincar\'e lemma in degree 
$q$.
\smallskip
By introducing a smooth Riemannian metric in $M$, 
and a smooth partition of unity 
$\{\chi_\ell,U_\ell\}$ subordinated to a coordinate 
trivialization atlas $(U_\ell,x_\ell)$ for $E^q$, 
we define for any compact subset $K$ of $M$
the seminorm $\| \,\cdot\,\|_{q,K,m}$ by\,:
\form{\dsize\|f\|_{q,K,m}=\sum_{\ell}\sup_{x\in K}\sum_{|\alpha|\leq m}
\left|\frac{\partial^{|\alpha|}(\chi_\ell f_\ell)(x)}{
\partial [x_{\ell}^1]^{\alpha_1}\cdots
\partial [x_{\ell}^N]^{\alpha_N}}\right|_{r_q}}
where $f_\ell(x)\in\Bbb{C}^{r_q}$ is the trivialization of
$f$ in $U_\ell$ and $|\,\cdot\,|_{r_q}$ is the standard Euclidean norm 
in $\Bbb{C}^{r_q}$. \par
Next we introduce smooth Hermitian scalar products $(\,\cdot\,|\,\cdot\,)_q$
on the fibers of the
$E^q$'s. This allows us to define the $L^2$-scalar product of
$f,g\in\Cal{E}^q(U)$ by\,:
\form{\int_U (f|g)_q d\lambda\, ,}
where $d\lambda$ is the measure associated to the Riemannian metric
of $M$. 
It is well defined when 
$\roman{supp}(f)\cap\roman{supp}(g)$ is compact in $U$.
Using these scalar products we obtain the 
(formal) {\it adjoint complex}
\form{\CD
\Cal{C}^{\infty}(M,E^0)@<{A_0^*}<<\Cal{C}^{\infty}(M,E^1)@<{A_1}<<
\Cal{C}^{\infty}(M,E^2)@<<<
\cdots\endCD}\edef\complessagg{\number\q.\number\t}
by defining the partial differential operator
$A_q^*:\Cal{C}^\infty(M,E^{q+1})@>>>\Cal{C}^{\infty}(M,E^q)$ by\,:
\form{\matrix\format\l\\
\int_M (A^*_qf|g)_q d\lambda\, = \,
\int_M (f|A_qg)_{q+1} d\lambda \\\vspace{2\jot}
\qquad
\forall f\in\Cal{C}^\infty(M,E^{q+1}),\;\forall
g\in \Cal{C}^\infty_{\roman{comp}}(M,E^{q+1})
\endmatrix}
where $\Cal{C}^\infty_{\roman{comp}}(M,E^{q+1})$ is the subspace
of $g\in \Cal{C}^\infty(M,E^{q+1})$ with $\roman{supp}(g)$ compact
in $M$.
%%%%%%%%%%%%%%%%%%%%%%%%%%%%%%%%%%%
\par In [AFN] the following was proved:
\thm{Assume that the complex
(\complesso)
admits the Poincar\'e lemma in degree $q$
at the point $x_0\in M$, i.e. we assume that
the sequence:
\form{\CD
\Cal E^{q-1}_{x_0}@>{A_{q-1}}>>\Cal E^q_{x_0}@>{A_{q}}>>\Cal E^{q+1}_{x_0}\, 
\endCD}
is exact.
Then, for every open neighborhood $\omega$ of $x_0$ in $M$ 
we can
find an open neighborhood $\omega_1$ of $x_0$ in $\omega$
such that 
\form{\forall f\in\Cal{E}^q(\omega)\;\roman{with}\;
A_qf=0\;\roman{in}\;\omega\quad\exists\; u\in\Cal{E}^{q-1}(\omega_1)\;
\roman{with}\; A_{q-1}u=f\;\roman{in}\;
\omega_1\,.
}\xdef\formulafai{\number\q.\number\t}}\edef\teorba{\number\q.\number\x}
As a consequence of the open mapping theorem for Fr\'echet spaces
we also have\,:
\thm{Let $\omega_1\subset\omega\subset M$ be open subsets 
such that the restriction map 
$\Cal{E}^*(\omega)@>>>\Cal{E}^*(\omega_1)$ induces the
zero map in cohomology\,:
\form{\{H^q(\omega,\Cal{E}^*,A_*) @>>>H^q(\omega_1,\Cal{E}^*,A_*)\}
=0\,,}
i.e. such that (\formulafai) holds true. Then,
for every compact $K_1\subset\omega_1$ and
every integer $m_1\geq 0$ there are a compact $K\subset\omega$,
an integer $m\geq 0$ and a constant $c>0$ such that\,:
the function $u$ in (\formulafai) can be chosen to satisfy\,:
\form{\|u\|_{q-1,K_1,m_1}\leq c\|f\|_{q,K,m}\, .}}
%%%%%%%%%%%%%%%%%%%%%%%%%%%%%%%%%%%%%%%%%%%%%%%%%%%%%%%%%%%%%%
From the preceding Theorem  (see again [AFN]) one obtains:
\thm{Assume that for (\complesso) we have (\formulafai).
Then for every compact subset $K_1$ of $\omega_1$ we can find
a compact subset $K$ of $\omega$,
an integer $m\geq 0$, and
a constant $c>0$
such that
for every $f\in\Cal E^q(\omega)$ with $A_qf=0$ in $\omega$ and
for every $v\in\Cal E^{q}(\omega_1)$ with
$\roman{supp}(v)\subset K_1$
we have:
\form{
\left|\dsize\int_{\omega_1}{(f|v)_q\, d\lambda}\right| \,\leq\,
c\,\|A^*_qv\|_{q-1,K_1,0}\, 
\|f\|_{q,K,m}\,. }} \edef\teoremafab{\number\q.\number\x}
Using the Hermitian inner product on the fibers of $E^q$ and 
the duality pairing associated to the $L^2$-product, we can define
for every open $U\subset M$ the space 
$\Cal{D'}^q(U)$ of $E^q$-valued distributions
in $U$. We note that to both (\complesso) and (\complessagg) 
we can associate 
complexes of partial differential operators on distributions,
which are compatible with the natural inclusion map
$\Cal{E}^q(U)\hookrightarrow\Cal{D'}^q(U)$ for
every $q=0,1,\hdots$ and every open $U\subset M$.\par
We denote by 
\form{H^q(U,\Cal{D'}^*,A_*)=\dsize\frac{\ker\left(A_q:\Cal{D'}^q(U)
@>>>\Cal{D'}^{q+1}(U)\right)}{\roman{Image}
\left(A_{q-1}:\Cal{D'}^{q-1}(U)
@>>>\Cal{D'}^{q}(U)\right)} }
and
\form{
H^q((x_0),\Cal{D'}^*,A_*)=\dsize\frac{\ker\left(A_q:\Cal{D'}^q_{x_0}
@>>>\Cal{D'}^{q+1}_{x_0}\right)}{\roman{Image}
\left(A_{q-1}:\Cal{D'}^{q-1}_{x_0}
@>>>\Cal{D'}^{q}_{x_0}\right)}
=\varinjlim_{U\ni x_0}H^q(U,\Cal{D'}^*,A_*)\,, 
}
the corresponding cohomology groups.
The natural inclusion maps 
$\Cal{E}^q(U)\hookrightarrow\Cal{D'}^q(U)$, 
induce natural maps in cohomology\,:
\form{H^q(U,\Cal{E}^*,A_*)@>>>H^q(U,\Cal{D'}^*,A_*)}
and
\form{ H^q((x_0),\Cal{E}^*,A_*)@>>>H^q((x_0),\Cal{D'}^*,A_*)
}\edef\coomindistrpunt{\number\q.\number\t}
for every open $U\subset M$ and $x_0\in M$.
\smallskip
Following the same arguments of [AFN] we can prove the following:
\thm{If (\coomindistrpunt), for fixed $q\geq 1$ and $x_0\in M$,
is the zero map, then for every open neighborhood $\omega$ of
$x_0$ in $M$ there exists an open neighborhood $\omega_1$ of
$x_0$ in $\omega$ such that\,:
\form{\forall f\in\Cal{E}^q(\omega)\;\roman{with}\;
A_qf=0\;\roman{in}\;\omega\quad\exists\; u\in\Cal{D'}^{q-1}(\omega_1)\;
\roman{with}\; A_{q-1}u=f\;\roman{in}\;
\omega_1\,.}}
Moreover, we have the analogue of Theorem \teoremafab\,:
\thm{Assume that for some $q\geq 1$ and 
open $\omega_1\subset\omega\subset M$ the composition\,:
\form{H^q(\omega,\Cal{E}^*,A_*)@>>>
H^q(\omega_1,\Cal{E}^*,A_*)@>>>H^q(\omega_1,\Cal{D'}^*,A_*)}
yields the zero map. Then for every compact subset $K_1\subset\omega_1$
there exist integers $m,m_1\geq 0$, a compact $K\subset\omega$
and a constant $c>0$ such that
for every $f\in\Cal E^q(\omega)$ with $A_qf=0$ in $\omega$ and
for every $v\in\Cal E^{q}(\omega_1)$ with
$\roman{supp}(v)\subset K_1$
we have\,:
\form{
\left|\dsize\int_{\omega_1}{(f|v)_q\, d\lambda}\right| \,\leq\,
c\,\|A^*_qv\|_{q-1,K_1,m_1}\, 
\|f\|_{q,K,m}\,. }} \edef\teoremafad{\number\q.\number\x}
%%%%%%%%%%%%%%%%%%%%%%%%%%%%%%%%%%%%%%%%%%%%%%%%%%%%%%%%%%%%%
\se{Preliminaries on $CR$ manifolds and notation}
In this paper $M$ will be a smooth ($\Cal C^\infty$) paracompact
manifold, of real dimension $2n+k$, with a smooth
$CR$ structure of type $(n,k)$: $n$ is its
{\it complex} $CR$ dimension and
$k$ its {\it real} $CR$ codimension. As an
abstract $CR$ manifold $M$ is a triple $\bold M=(M,HM,J)$,
where $HM$ is a smooth real vector subbundle of rank $2n$ of the
real tangent bundle $TM$, and where $J:HM@>>>HM$ is a smooth fiber
preserving isomorphism such that $J^2=-I$. It is also required that
the {\it formal integrability conditions}
 $\left[\Cal{C}^{\infty}(M,T^{0,1}M),\Cal{C}^{\infty}(M,T^{0,1}M)\right]
\subset\Cal{C}^{\infty}(M,T^{0,1}M)$
be satisfied. Here
 $T^{0,1}M=\left\{X+iJX\; \right|\; \left. X\in HM\right\}$ is the
complex subbundle of the complexification $\Bbb C HM$ of $HM$
corresponding to the eigenvalue $-i$ of $J$; we have 
 $T^{1,0}M\cap T^{0,1}M= 0$ and $T^{1,0}M\oplus T^{0,1}M=\Bbb C HM$,
where $T^{1,0}M=\overline{T^{0,1}M}$. When $k=0$, we recover the
abstract definition of a {\it complex} manifold, via the 
Newlander-Nirenberg theorem.\par
We denote by
 $H^0M=\{\xi\in T^*M\, | \, \langle X,\xi\rangle=0\quad\forall X\in 
H_{\pi(\xi)}M\}$ 
the {\it characteristic bundle} of $M$. To each $\xi\in H_x^0M$,
we associate the Levi form at $\xi$:
\form{\Cal{L}_{\xi}(X)=\xi([J\tilde X,\tilde X])=d\tilde\xi(X,JX) \quad
\text{for}\quad X\in H_xM\,}
which is Hermitian for the complex structure of $H_xM$ defined by $J$.
Here $\tilde\xi$ is a section of $H^0M$ extending $\xi$ and
 $\tilde X$ a section of $HM$ extending $X$.\par
Let $\Cal E^*(M)=\bigoplus_{h=0}^{2n+k}{\Cal E^{(h)}(M)}$ denote the
Grassmann algebra of smooth, complex valued differential forms on $M$.
We denote by $\Cal J$ the ideal of $\Cal E^*(M)$ that annihilates
 $T^{0,1}M$:
\form{\Cal J=\left.\left\{
\alpha\in\bigoplus_{h\geq 1}{\Cal E^{(h)}(M)}\; \right|\;
\alpha|_{T^{0,1}M}=0\,\right\}\, .}
By the formal integrability conditions we have $d\Cal J\subset\Cal J$.
We also consider the powers $\Cal J^p$ of the ideal $\Cal J$,
obtaining a decreasing sequence of $d$-closed ideals of $\Cal E^*(M)$:
\form{\Cal E^*(M)=\Cal J^0\supset\Cal J^1\subset\Cal J^2\supset\cdots
\supset\Cal J^{n+k-1}\supset\Cal J^{n+k}\supset\Cal J^{n+k+1}=\{0\}\, .}
Passing to the quotients, the exterior differential defines linear maps:
\par\noindent
 $\bar\partial_M:\Cal J^p/\Cal J^{p+1}@>>>\Cal J^p/\Cal J^{p+1}$.
The grading of $\Cal E^*(M)$ induces a grading of $\Cal J^p/\Cal J^{p+1}$:
\form{\Cal J^p/\Cal J^{p+1}=\bigoplus_{j=0}^n{\Cal Q_M^{p,j}(M)}\, .}
As $\bar\partial_M\circ\bar\partial_M=0$, we obtain the 
{\it tangential Cauchy-Riemann complexes} for $0\leq p\leq n+k$:
\form{(\Cal Q^{p,*},\bar\partial_M)=\left\{
0@>>>\Cal Q^{p,0}_M(M)@>{\bar\partial_M}>>
\Cal Q^{p,1}_M(M)@>>> \cdots @>{\bar\partial_M}>>
\Cal Q^{p,n}_M(M)@>>>0\right\}\, .}\edef\dcomplesso{\number\q.\number\t}
We also note that $\Cal Q^{p,j}_M(M)=\Cal{C}^{\infty}(M,Q_M^{p,j})$ for
complex vector bundles $Q_M^{p,j}$ on $M$ of rank
  $\binom{n+k}{p}+\binom{n}{j}$ and (\dcomplesso) is a complex of partial
differential operators of the first order. We denote by
$\Omega^p_M$ the sheaf of germs of smooth sections $f$ of $Q^{p,0}_M$
satisfying $\bar\partial_Mf=0$. Note that $Q^{0,0}_M$ is the trivial
complex line bundle over $M$; we set $\Omega^0_M=\Cal O_M$ and call it
the sheaf of germs of smooth $CR$ functions on $M$.
\par
If $(M,HM,J_M)$ and $(N,HN,J_N)$ are two $CR$ manifolds, we say that
a differentiable map $\phi:M@>>>N$ is $CR$ iff: 
($i$) $d\phi(HM)\subset HN$; ($ii$) $d\phi(J_MX)=J_Nd\phi(X)$ for all
 $X\in HM$. \par
We say that a $CR$ manifold 
$(M,HM,J_M)$ 
of type $(n,k)$ is {\it locally embeddable} if for each point
 $x\in M$ we can find an open neighborhood $U$ of $x$ in $M$, 
an open subset $\tilde U$ of $\Bbb C^{n+k}$, and a smooth $CR$
map $\phi:U@>>>\tilde U$ which is an embedding.
\medskip
For a locally embeddable $M$ we shall give now a 
description in local coordinates of the Levi form and of the
tangential Cauchy-Riemann complex. \par
Let $x\in M$ and $U$ open in $\Bbb C^{n+k}$ be as above. We can assume
that $x$ is the origin $0$ of $\Bbb C^{n+k}$ and that
\form{M\cap U=\{x\in U\, | \, \rho_1(z)=0,\,\hdots ,\,\rho_k(z)=0\}}
where $\rho_1,\,\hdots ,\,\rho_k$ are real valued smooth functions 
on $U$ and 
\form{\partial\rho_1(z)\wedge\cdots\wedge\partial\rho_k(z)\neq 0\qquad
\text{for}\quad z\in U\, .}
The holomorphic tangent space to $M$ at a point $z\in M$, having chosen
holomorphic coordinates on $T_z\Bbb{C}^{n+k}$, is identified with
$T^{1,0}_zM$ and
is described by:
\form{H_zM\,
\simeq T^{1,0}_zM\,=\,\left\{u=(u^\alpha)\in\Bbb C^{n+k}\, \left| \, 
\dsize\sum_{\alpha=1}^{n+k}{\frac{\partial\rho_j(z)}{\partial z^\alpha}
u^{\alpha}}=0\;\text{for}\; j=1,\hdots,k\right\}\right. .}
We also have:
\form{H_z^0M\,=\,\left.\left\{{\xi=\left(
\sum_{j=1}^k{\lambda^jd^c\rho_j(z)}\right)}
\, \right|\, \lambda^1,\hdots ,\lambda^k\in\Bbb R\right\}}
and the Levi form at $\xi=\left(\sum_{j=1}^k{\lambda^jd^c\rho_j(z)}\right)$
is the complex Hessian
\form{\dsize\sum_{j=1}^k
\dsize\sum_{\alpha,\beta=1}^{n+k}{\lambda^j\frac{\partial\rho_j(z)}{
\partial z^\alpha\partial\bar z^\beta}u^\alpha\bar u^\beta}\qquad
\text{for}\quad u=(u^\alpha)\in T^{1,0}_zM\simeq H_zM\, .}
\smallskip
Let $z^{n+j}=t^j+is^j$ with $t^j,s^j\in\Bbb R$, for $j=1,\hdots,k$.
By a linear change of coordinates we can obtain that near $x=0$:
\form{\rho_j(z)=s^j-h_j(z^1,\hdots,z^n,t^1,\hdots,t^k),\;
\text{with}\;h_j=0(2)\;\text{for}\; j=1,\hdots, k\, .}
Near $x=0$, the complex coordinates $z^1,\hdots,z^n$ and the
real coordinates $t^1,\hdots, t^k$ define smooth coordinates on $M$.
\par
Denote by $C=(C_{\alpha,j})$ the matrix
\form{C=-2i\left(I+i\left(\frac{\partial 
h_\ell(z)}{\partial t^j}\right)\right)^{-1}\left(\frac{\partial h_j(z)}{
\partial\bar z^\alpha}\right)}
and consider the partial differential operators
\form{\bar L_\alpha=\frac{\partial}{\partial\bar z^\alpha}+
\dsize\sum_{j=1}^k{C_{\alpha}^j\frac{\partial}{\partial t^j}}\,.}
On a neighborhood $V$ of $0$ in $M$ the natural pull-back composed with the
projection onto the quotient define an identification:
\form{\Cal Q^{p,q}_M(V)\simeq\left\{
\dsize\sum_{\smallmatrix
1\leq \alpha_1<\cdots < \alpha_p\leq n+k\\
1\leq \beta_1\leq\cdots < \beta_q\leq n
\endsmallmatrix}{\!\!\!\!\!\!\!\!\!\!\!\!
a_{\alpha_1,\hdots,\alpha_p;\beta_1\hdots\beta_q}
dz^{\alpha_1}\wedge\cdots\wedge dz^{\alpha_p}
\wedge d\bar z^{\beta_1}\wedge\cdots\wedge d\bar z^{\beta_q}}\right\}}
where $a_{\alpha_1,\hdots,\alpha_p;\beta_1\hdots\beta_q}=
a_{\alpha_1,\hdots,\alpha_p;\beta_1\hdots\beta_q}
(z^1,\hdots,z^n,t^1,\hdots,t^k)$ are smooth complex valued functions
on $V$.
\par
With this identification, and observing that $[\bar L_\alpha,\bar L_\beta]=0$
for $1\leq\alpha,\beta\leq n$, we obtain:
\form{\matrix\format\l\\
\bar\partial_M\left(\dsize\sum_{\smallmatrix
1\leq \alpha_1<\cdots < \alpha_p\leq n+k\\
1\leq \beta_1\leq\cdots < \beta_q\leq n
\endsmallmatrix}{a_{\alpha_1,\hdots,\alpha_p;\beta_1\hdots\beta_q}
dz^{\alpha_1}\wedge\cdots\wedge dz^{\alpha_p}
\wedge d\bar z^{\beta_1}\wedge\cdots\wedge d\bar z^{\beta_q}}\right)\\
{}\\
%=
{}\\
%\dsize\sum_
\underset{\smallmatrix
1\leq \alpha_1<\cdots < \alpha_p\leq n+k\\
1\leq \beta_1\leq\cdots < \beta_q\leq n\\
1\leq\beta\leq n
\endsmallmatrix}\to{=\quad\dsize\sum}{\left(\bar L_\beta
a_{\alpha_1,\hdots,\alpha_p;\beta_1\hdots\beta_q}\right)
d\bar z^\beta\wedge dz^{\alpha_1}\wedge\cdots\wedge dz^{\alpha_p}
\wedge d\bar z^{\beta_1}\wedge\cdots\wedge d\bar z^{\beta_q}}\, .
\endmatrix}

%%%%%%%%%%%%%%%%%%%%%%%%%%%%%%%%%%%%%%%%%%%%%%%%%%%%%%%%%%%%%%%%
\se{The $CR$ structure of the characteristic bundle}
In this section we define {\it regular points} of
$H^0M$ and
prove a result that was formulated in \cite{Tu,
Theorem 2} for a generically embedded
$CR$ manifold. In that paper we could not find a proof of
the stated 
Levi flatness of the $CR$ structure of the conormal bundle.
This fact is interesting for our discussion,
because it
entails the {\it regularity} of the characteristic
codirections $\xi\in H^0M$ at which the 
rank of the Levi form is maximal.
This regularity is
requested in Theorems 5.1 and 5.2
to discuss
the non validity of the 
Poincar\'e Lemma. \par
%%%%%%%%%%%%%%%%%%%%%%%%%%%%%%%%%%%%%%%%%%%%%%%%%%%%%%%
We have preferred to consider the
characteristic bundle, which is more intrinsically related to
the differential geometry of (abstract) $CR$ manifolds,
rather than the conormal bundle of \cite{Tu}.
However, our results imply those of \cite{Tu, Theorem 2}
since, when $M$ is a
generic $CR$ submanifold of a complex manifold $\Tilde{M}$, 
the characteristic bundle $H^0M$ is the image, by the dual map
$\Cal{J}^*$ of the complex structure $\Cal{J}:T\Tilde{M}@>>>
T\Tilde{M}$, of the conormal bundle $\Tilde{M}^*(M)$ of
$M$ in $\Tilde{M}$\,:
$H^0M=\Cal{J}^*\left(\Tilde{M}^*(M)\right)$, and
$\Cal{J}^*:T^*\Tilde{M}@>>>T^*\Tilde{M}$ is biholomorphic.
\smallskip
Let $M$ be smooth abstract 
$CR$ submanifold, of type $(n,k)$. Denote by
$\vartheta$ the tautological $1$-form on $H^0M$. If
$\pi:H^0M @>>> M$ is the natural projection, we have\,:
\form{\vartheta(X)=\xi(\pi_*(X))\quad\text{if}\quad \xi\in H^0_xM
\quad\text{and}\quad X\in T_{\xi}H^0M\,.}
For each point $\xi\in H^0M$, we set\,:
\form{\cases
N_\xi  H^0M\,=\,\{X\in T_{\xi}H^0M\, |\, \vartheta(X)=0\,\}\,,\\
H_\xi H^0M\,=\,\{X\in N_{\xi}H^0M\, |\, 
d\vartheta(X,Y)=0\quad\forall Y\in T_\xi  H^0M\,\}\,.
\endcases}
Let $\Dot{H}^0M$ be the open subset of the nonzero cotangent
vectors in ${H}^0M$.
Then $\xi@>>>N_\xi  H^0M$ is a smooth distribution of hyperplanes
in $T\dot{H}^0M$.  We note however, as the discussion below will
clarify,
that in general 
the dimension of $H_\xi H^0M$ is not constant and
therefore $\xi @>>>H_\xi H^0M$ may fail  to be a smooth
subbundle of $T\dot{H}^0M$.\par  
If $U$ is an open subset of $H^0M$, we define\,:
\form{\cases
\frak{N}(U)=\{X\in\Cal{C}^\infty(U,TH^0M)\, | \,
X(\xi)\in N_\xi H^0M \quad\forall\xi\in U\,\}\, ,\\
\frak{H}(U)=\{X\in\Cal{C}^\infty(U,TH^0M)\, | \,
X(\xi)\in H_\xi H^0M \quad\forall\xi\in U\,\}.
\endcases}
\lem{Let $U$ be an open subset of $\Dot{H}^0M$. Then\,:
\form{
\cases
[\frak{H}(U),\frak{N}(U)]\subset \frak{N}(U)\,,\\
[\frak{H}(U),\frak{H}(U)]\subset \frak{H}(U).
\endcases}
}\edef\lemmaa{\number\q.\number\x}
\dimo
For $X\in\frak{H}(U)$ and $Y\in\frak{N}(U)$, we have\,:
$$\vartheta([X,Y])=-d\vartheta(X,Y)+X\vartheta(Y)-Y\vartheta(X)=0\,.$$ 
Indeed $d\vartheta(X,Y)=0$ because $X(\xi)\in H_{\xi}H^0M$, and
$\vartheta(Y)=0$, $\vartheta(X)=0$ because 
$Y(\xi),\;X(\xi)\in N_{\xi}H^0M$, for all $\xi\in U$.
This shows that $[\frak{H}(U),\frak{N}(U)]\subset \frak{N}(U)$.
\par
Let now $X,Y\in\frak{H}(U)$. Then
$[X,Y]\in \frak{N}(U)$ and, to complete the proof of the Lemma,
we need only to verify that
$d\vartheta([X,Y],Z)=0$ for all $Z\in\Cal{C}^\infty(U,TH^0M)$.
Since $\vartheta\neq 0$ at each point $\xi\in U$, 
we can as well assume that $\vartheta(Z)$ is constant in $U$.
We have\,:
$$\matrix\format\r&\l\\
0=d\vartheta(X,Z)&=X\vartheta(Z)-Z\vartheta(X)-\vartheta([X,Z])\\
&=-\vartheta([X,Z])
\endmatrix$$
and likewise $\vartheta([Y,Z])=0$.
Therefore
$$\vartheta([X,[Y,Z]])=-d\vartheta(X,[Y,Z])+X\vartheta([Y,Z])-
[Y,Z]\vartheta(X)=0$$
and likewise
$\vartheta([Y,[X,Z]])=0$.
Hence we have\,:
$$\matrix\format\r&\l\\
d\vartheta([X,Y],Z)&=[X,Y]\vartheta(Z)-Z\vartheta([X,Y])-
\vartheta([[X,Y],Z])\\\vspace{2\jot}
&=-\vartheta([X,[Y,Z]])+\vartheta([Y,[X,Z]])=0\, .
\endmatrix
$$
This shows that $[\frak{H}(U),\frak{H}(U)]\subset \frak{H}(U)$.
\qed
\medskip
Assume now that that we are given a generic $CR$ embedding
$M\hookrightarrow \Tilde{M}$ of $M$ into an $(n+k)$-dimensional complex
manifold $\Tilde{M}$. This yields an embedding 
$H^0M\hookrightarrow T^*\Tilde{M}$. The cotangent bundle $T^*\Tilde{M}$
has a complex structure, that we denote by 
$\frak{J}\,: TT^*\Tilde{M}@>>> TT^*\Tilde{M}$. 
We shall still denote by $\pi:T^*\Tilde{M}@>>>\Tilde{M}$ the natural
projection, extending $H^0M@>\pi>>M$.
\lem{Let $x\in M$ and $\xi\in H^0_xM\setminus\{0\}$. Then
\form{H_{\xi}H^0M=T_{\xi}H^0M\cap\frak{J}T_{\xi}H^0M\, .}
}\edef\lemmab{\number\q.\number\x}
\dimo
In the proof we can as well assume that $\Tilde{M}$ is an open
subset of $\Bbb{C}^{n+k}$. In particular we shall utilize the
standard trivialization of the cotangent bundle and identify 
$H^0M$ with a submanifold of 
$\Tilde{M}\times\left(\Bbb{R}^{2n+2k}\right)^*$\,:
\form{H^0M\, =\, \left.\left\{
(x,x^*)\in M\times \left(\Bbb{R}^{2n+2k}\right)^*\, \right|\,
\langle x^*,H_xM\rangle=\{0\}\,\right\}\, .}
We can assume that
$M$  is defined in $\Tilde{M}$ by
\form{M=\{x\in \Tilde{M}\, | \, \rho_1(x)=0,\;\hdots,\; \rho_k(x)=0\}\,
,}\edef\defemme{\number\q.\number\t}
with $\rho_1,\;\hdots,\;\rho_k\in\Cal{C}^\infty(\Tilde{M},\Bbb{R})$ and
$\partial\rho_1(x)\wedge\cdots\wedge\partial\rho_k(x)\neq 0$ for all
$x\in \Tilde{M}$. \par
With $d^c=(\partial-\bar\partial)/i$, we obtain\,:
%%%%%%%%%%%%%%%%%%%%%%%%%%%%%%%%%%%%%%%%%%%%%%%%%%%%%%%%%%%%%
\form{H^0M\, = \,
\left.\left\{\left(x,\sum_{j=1}^k{a^jd^c\rho_j(x)}\right)\,\right|\,
x\in M\,,\;a^1,\hdots,a^k\in\Bbb{R}\right\}\,.}
In this way a tangent vector to $H^0M$ at 
$\xi_0=(x_0,{x^*}^0)=\left(x_0,\sum_{j=1}^k{a^j_0d^c\rho_j(x_0)}\right)$ 
is 
identified to a vector in 
$\Bbb{R}^{2n+2k}\times\left(\Bbb{R}^{2n+2k}\right)^*$
of the form\,:
\form{X=\left(v_0,\sum_{j=1}^k{\lambda^j \langle d^c\rho_j(x_0),
\,\cdot\,\rangle}+
\sum_{j=1}^k{a^j_0 dd^c\rho_j(x_0) (v_0,\,\cdot\,)}\right)
\, ,}\edef\formulag{\number\q.\number\t}
with $v_0\in T_{x_0}M$ and $\lambda^1,\hdots,\lambda^k\in\Bbb{R}$. 
We recall the notation
$\Cal{J}$  for the complex structure of $\Tilde{M}$, and hence here of
$\Bbb{C}^{n+k}$, and
$\frak{J}$ for the complex structure of $T^*\Tilde{M}$, and hence here of
$T^*\Bbb{C}^{n+k}$.
We have\,:
\form{\frak{J}X=
\left(\Cal{J}v_0,
\sum_{j=1}^k{\lambda^j \langle d^c\rho_j(x_0),\,\Cal{J}\,\cdot\,
\rangle}+
\sum_{j=1}^k{a^j_0 dd^c\rho_j(x_0)(v_0,\,
\Cal{J}\,\cdot\,)}\right)\, .}
Therefore $X$ and $\frak{J}X$ both belong
to $T_{\xi}H^0M$ if, and only if, 
$$v_0\in H_{x_0}M \, ,\tag $i$ $$
and there exist real numbers $\mu^1,\hdots,\mu^k$ such that\,:
$$ \matrix\\
\sum_{j=1}^k{\lambda^j \langle d^c\rho_j(x_0),\Cal{J}w\,\rangle}+
\sum_{j=1}^k{a^j_0 dd^c\rho_j(x_0)(v_0,\,\Cal{J}w\,)}\\\vspace{2\jot}
\qquad\qquad =
 \sum_{j=1}^k{\mu^j \langle d^c\rho_j(x_0),\,w\,\rangle}+
\sum_{j=1}^k{a^j_0 dd^c\rho_j(x_0)(\Cal{J}v_0,\,w\,)}\\\vspace{2\jot}
\qquad\qquad\qquad 
\forall w\in\Bbb{R}^{2n+2k}\endmatrix
\tag $ii$ 
$$
From this we deduce that\,:
\form{\matrix
\sum_{j=1}^k{a^j_0\left(dd^c\rho_j(x_0)(v_0,\Cal{J}w)- 
dd^c\rho_j(x_0)(\Cal{J}v_0,w)\right)}\\
\qquad\qquad =\sum_{j=1}^k{\mu^j\,\langle
d^c\rho_j(x_0),w\rangle}\qquad\forall w\in T_{x_0}M\,
\endmatrix}\edef\formulai{\number\q.\number\t}
and 
\form{\matrix\format\l\\
\sum_{j=1}^k{a^j_0\left(dd^c\rho_j(x_0)(v_0,\Cal{J}w)- 
dd^c\rho_j(x_0)(\Cal{J}v_0,w)\right)}
\\
\qquad =-\sum_{j=1}^k{\lambda^j\,\langle
d^c\rho_j(x_0),\Cal{J}w\rangle}
=-\sum_{j=1}^k{\lambda^j\,\langle
d\rho_j(x_0),w\rangle}\\
\qquad\qquad\qquad\qquad\qquad\qquad\forall w\in 
\left(T_{x_0}M\right)^\perp\,
\endmatrix
}\edef\formulaj{\number\q.\number\t}
Here the perpendicular is taken with respect to the standard
metric of $\Bbb{R}^{2n+2k}$. \par
We can find the real numbers
$\mu^j$ to satisfy (\formulai) if and only if 
$$
\sum_{j=1}^k{a^j_0\left(dd^c\rho_j(x_0)(v_0,\Cal{J}w)- 
dd^c\rho_j(x_0)(\Cal{J}v_0,w)\right)}=0\qquad\forall w\in H_{x_0}M\,.
$$
When $v_0,w\in H_{x_0}M$ we have 
$$dd^c\rho_j(x_0)(v_0,\Cal{J}w)=- 
dd^c\rho_j(x_0)(\Cal{J}v_0,w)\qquad\forall j=1,\hdots,k
$$
and therefore (\formulai) is equivalent to
\form{v_0\in \ker\Cal{L}_{\xi_0}\,.}
Finally we note that (\formulaj) uniquely determines the 
real coefficients
$\lambda^j$ for any given $v_0\in\ker\Cal{L}_{\xi_0}$.
\par
In this way we characterized a vector 
$X\in T_{\xi_0}H^0M\cap \frak{J}T_{\xi_0}H^0M$:
it has the form (\formulag) with
\form{2\sum_{j=1}^k{a^j_0dd^c\rho_j(v_0,\Cal{J}w)}+\sum_{j=1}^k{\lambda^j
\langle d\rho_j(x_0),w\rangle }=0\qquad\forall w\in\left(T_{x_0}M\right)^\perp
\, . }\edef\caratterizzazione{\number\q.\number\t}
\smallskip
Next we compute $H_{\xi_0}H^0M$. To this aim we observe that, denoting
by $x^*_\alpha$ the dual real variables of $\left(\Bbb{R}^{2n+2k}\right)^*$
we can write a tangent vector in the form\,:
\form{X=\sum_{\alpha=1}^{2n+2k}{v^\alpha\frac{\partial}{\partial x^\alpha}}
+\sum_{j=1}^k{\lambda^j\langle d^c\rho_j(x),e_\alpha\rangle 
\frac{\partial}{\partial x^*_\alpha}}+
\sum_{j=1}^k{a^j\,dd^c\rho_j(x)(v,e_\alpha) 
\frac{\partial}{\partial x^*_\alpha}}}
where $e_\alpha$ ($\alpha=1,\hdots,2n+2k$) is the canonical basis
of $\Bbb{R}^{2n+2k}\simeq\Bbb{C}^{n+k}$ and the vector
$v=\sum_{\alpha=1}^{2n+2k}v^\alpha(x,x^*)\frac{\partial}{\partial x^\alpha}$
belongs to  
$T_xM$. It is convenient to consider the coefficients
$v^\alpha$, $\lambda^j$ to be
smooth functions of $x,x^*$, defined everywhere in 
$\Dot{T}^*\Tilde{M}$. A vector field of this form belongs
to $\frak{N}(U)$ (for a neighborhood $U$ of $\xi_0$ in $H^0M$)
iff 
\form{\sum_{j=1}^k{a^j\langle d^c\rho_j(x),v\rangle}=0\quad \text{in}
 \quad U\,.}
By differentiation we obtain\,:
\form{\sum_{\alpha=1}^{2n+2k}\left(
\sum_{j=1}^k \frac{\partial a^j}{\partial x^*_\alpha}
\langle d^c\rho_j,v\rangle\, + \, \sum_{j=1}^k a^j\left\langle
d^c\rho_j,\frac{\partial v}{\partial x^*_\alpha}\right\rangle
\right)dx^*_\alpha=0\quad\text{in}\quad U\, 
.}
We note that 
\form{x^*_\alpha=\sum_{j=1}^k{a^j\langle d^c\rho_j(x),e_\alpha\rangle}
\quad\text{in}\quad U}
and therefore
\form{dx^*_\alpha=
\sum_{\beta=1}^{2n+2k}\sum_{j=1}^k{\frac{\partial a^j}{\partial
x^*_\beta}\langle d^c\rho_j,e_\alpha\rangle}
dx^*_\beta\,\quad\text{in}\quad U\, .}
Let
$$Y=\sum_{\alpha=1}^{2n+2k}{w^\alpha\frac{\partial}{\partial x^\alpha}}
+\sum_{j=1}^k{\mu^j\langle d^c\rho_j(x),e_\alpha\rangle 
\frac{\partial}{\partial x^*_\alpha}}+
\sum_{j=1}^k{a^j\,dd^c\rho_j(x)(w,e_\alpha) 
\frac{\partial}{\partial x^*_\alpha}} $$
be another vector field in $U$. 
Let $[X,Y]=\sum_{\alpha=1}^{2n+2k}{
\left(A^\alpha\frac{\partial}{\partial x^\alpha}+
B_\alpha \frac{\partial}{\partial x^*_\alpha}\right)}$.
Then\,: 
\form{\matrix\format\r&\l\\
A^\alpha\; =&\;\dsize
v^\beta\frac{\partial w^\alpha}{\partial x_\beta}-
w^\beta\frac{\partial v^\alpha}{\partial x_\beta} \\\vspace{2\jot}
&+\dsize\sum_{j=1}^k{\lambda^j\langle d^c\rho_j,e_\beta\rangle
\frac{\partial w^\alpha}{\partial x^*_\beta}}
-\sum_{j=1}^k{\mu^j\langle d^c\rho_j,e_\beta\rangle
\frac{\partial v^\alpha}{\partial x^*_\beta}}\\\vspace{2\jot}
&+\dsize\sum_{j=1}^k {a^j dd^c\rho_j(v,e_\beta)
\frac{\partial w^\alpha}{\partial x^*_\beta}}
-
\sum_{j=1}^k{ a^j dd^c\rho_j(w,e_\beta)
\frac{\partial v^\alpha}{\partial x^*_\beta}}\, .
\endmatrix}
Assume that $X\in\frak{N}(U)$ and impose the condition 
that $d\vartheta(X,Y)=0$ for all $Y\in\Cal{C}^\infty(U,TH^0M)$. 
We can choose $Y$ is such a way that $\vartheta(Y)=\text{constant}$
on a neighborhood of $U$ in $T^*\Tilde{M}$, so that\,:
 \form{\sum_{\alpha=1}^{2n+2k}\left(
\sum_{j=1}^k \frac{\partial a^j}{\partial x^*_\alpha}
\langle d^c\rho_j,w\rangle\, + \, \sum_{j=1}^k a^j\left\langle
d^c\rho_j,\frac{\partial w}{\partial x^*_\alpha}\right\rangle
\right)dx^*_\alpha=0\quad\text{in}\quad U\, 
.}
%%%%%%%%%%%%%%%%%%%% 
Then $d\vartheta(X,Y)=-\vartheta([X,Y])$ and we are lead to the
condition\,:
\form{\matrix\format\l\\
\sum_{j=1}^ka^j\langle d^c\rho_j,[v,w]\rangle\\\vspace{2\jot}
\quad
-\sum_{j=1}^k\lambda^j\langle d^c\rho_j,w'\rangle
+\sum_{j=1}^k\mu^j\langle d^c\rho_j,v'\rangle\\\vspace{2\jot}
\quad
-\sum_{j=1}^ka^j dd^c\rho_j(v,w')+\sum_{j=1}^k a^j dd^c\rho_j(w,v')
\\\vspace{3\jot}
=
-\sum_{j=1}^ka^j dd^c\rho_j(v,w)\\\vspace{2\jot}
\quad
-\sum_{j=1}^k\lambda^j\langle d^c\rho_j,w'\rangle
+\sum_{j=1}^k\mu^j\langle d^c\rho_j,v'\rangle\\\vspace{2\jot}
\quad
-\sum_{j=1}^ka^j dd^c\rho_j(v,w')+\sum_{j=1}^k a^j dd^c\rho_j(w,v')
=0\, ,
\endmatrix}\edef\formulas{\number\q.\number\t}
where
$$v'=\sum_{j=1}^k\sum_{\alpha=1}^{2n+2k}
\frac{\partial a^j}{\partial x^*_\alpha}\langle d^c\rho_j,v\rangle e_\alpha\,
,
\qquad
w'=\sum_{j=1}^k\sum_{\alpha=1}^{2n+2k}
\frac{\partial a^j}{\partial x^*_\alpha}\langle d^c\rho_j,w\rangle e_\alpha
$$
are the components of $v$ and $w$ orthogonal to $HM$. The vector
$X_{\xi_0}$ belongs to $H_{\xi_0}H^0M$ if and only if
(\formulas) holds at the point $\xi_0=(x_0,{x^*}^0)$ for all choices
of $\mu^1_0,\hdots,\mu^k_0\in\Bbb{R}$ and for every $w_0\in T_{x_0}M$.
Taking first $w_0=0$ and arbitrary $\mu^i_0$, we obtain\,:
$$\langle d^c\rho_j(x_0),v'(\xi_0)\rangle=0\quad
\text{for}\quad j=1,\hdots,k\,,$$
i.e. $v(\xi_0)\in H_{x_0}M$. Next, letting $w_0$ vary in 
$H_{x_0}M$, we obtain that\,:
$$\sum_{j=1}^ka^j_0 dd^c\rho_j(v({\xi_0}),w_{0})=0\qquad
\forall w_{0}\in H_{x_0}M$$
and hence $v({\xi_0})\in \ker\Cal{L}_{\xi_0}$.
Finally,
we have\,:
$$2\sum_{j=1}^ka^j_0dd^c\rho_j(v(\xi_0),w'_0)+\sum_{j=1}^k\lambda^j_0
\langle d^c\rho_j(x_0),w'_0\rangle=0 \qquad\forall w'_0\in T_{x_0}M\cap
\left(H_{x_0}M\right)^\perp\, .$$
Since $\langle d^c\rho_j,\Cal{J}w\rangle=
\langle d\rho_j,w\rangle$, and $\Cal{J}\left(\left(T_{x_0}M\right)^\perp
\right)=T_{x_0}M\cap
\left(H_{x_0}M\right)^\perp$, by (\caratterizzazione)
the proof of the lemma is complete. \qed
\medskip
\thm{Let $M$ be a generic $CR$ submanifold, of type $(n,k)$. 
Assume that the Levi form
$\Cal{L}_{\xi}$ has constant rank $m$ for all $\xi$ in an open subset
$U$ of $\Dot{H}^0M$.  Then $U$ is a $CR$ submanifold of
$T^*\Tilde{M}$, of type $(n-m,2k+2m)$, which is Levi flat and
foliated by complex submanifolds of dimension $(n-m)$.
Note that the embedding $U\hookrightarrow T^*\Tilde{M}$ is
not generic when $m<n$.}
\dimo We have\,: ($i$) $[\frak{N}(U),\frak{N}(U)]
\subset\frak{N}(U)$
by Lemma \lemmaa; ($ii$) $\frak{J}(\frak{N}(U))=\frak{N}(U)$,
by Lemma \lemmab; 
($iii$) $\frak{N}(U)$ is a distribution
of constant rank $(2n-2m)$ by Lemma {\lemmab} . Hence 
the real Frobenious theorem
provides the foliation, and by the classical theorem of
Levi-Civita, each leaf is a complex 
submanifold of dimension $(n-m)$. \qed
\medskip\edef\teoremaa{\number\q.\number\x}
We say that a point $\xi_0\in\Dot{H}^0M$ is {\it regular} if
there exists an open neighborhood $U$ of $\xi_0$ in $\Dot{H}^0M$
and a smooth submanifold $V$ of $U$ with
\form{\xi_0\in V\subset U\,,\quad  T_{\xi}V\subset H_{\xi}H^0M\;\forall
\xi\in V\, ,
\;\quad\text{and}\quad
 T_{\xi_0}V=H_{\xi_0}H^0M\, .}
Since the rank of the Levi form
$\Cal{L}_{\xi}$ is a lower
semicontinuous function of $\xi\in\Dot{H}^0M$,
Theorem {\teoremaa } yields\,:
\cor{Let $M$ be a generic $CR$ submanifold, of type $(n,k)$, 
and $\xi_0\in\Dot{H}^0M$ with
$\roman{rank}\Cal{L}_{\xi_0}=\limsup_{\xi@>>>\xi_0}
\roman{rank}\Cal{L}_{\xi}$.
Then $\xi_0$ is regular and there is an open neighborhood
$U$ of $\xi_0$ in $\Dot{H}^0M$ consisting of regular points
$\xi$ with
$\roman{rank}\Cal{L}_{\xi}=\roman{rank}\Cal{L}_{\xi_0}$.}
We have\,:
\lem{Let $M$ be a generic $CR$ submanifold, of type $(n,k)$, of a complex
manifold $\Tilde{M}$.
Let $\xi_0\in\Dot{H}^0M$ be regular. If the Levi form
$\Cal{L}_{\xi_0}$ has rank $m$, then there exists an $(n-m)$-dimensional
smooth complex submanifold $V^*$ of 
a neighborhood of $\xi_0$ in $T^*\Tilde{M}$, with\,:
\form{\xi_0\in V^*\subset H^0M\, .}
The projection $W=\pi(V^*)$ is a smooth complex submanifold 
of a neighborhood of $\pi(\xi_0)$ in
$\Tilde{M}$, contained in $M$.}
\dimo 
The dimension of $H_\xi H^0M$ is an upper semicontinuous function
of $\xi\in\Dot{H}^0M$. Therefore it 
remains constant and equal to $2(n-m)$
on an open neighborhood
of $\xi_0$ in $V$. Let $V^*$ be such a neighborhood. We have
$T_{\xi}V^*=H_{\xi}H^0M$ for all $\xi\in V^*$. Hence 
$\frak{J}T_{\xi}V^*=T_{\xi}V^*$ for all $\xi\in V^*$ and then,
by the theorem of Levi-Civita, $V^*$ is a complex submanifold
of dimension $(n-m)$ of a neighborhood of $\xi_0$ in
$T^*\Tilde{M}$. Since the fibers of $H^0M @>>>M$ are totally real,
$V^*$ is transversal to the fibers and therefore the map
$V^*\ni\xi @>>>\pi(\xi)\in M$ is a local diffeomorphism. 
It becomes a diffeomorphism after substituting to $V^*$ its intersection
with a suitable small neighborhood of $\xi_0$ in $T^*\Tilde{M}$.
Finally, $\pi(V^*)$ is a complex submanifold of an open neighborhood
of $\pi(\xi_0)$, contained
in $M$, because the
projection $T^*\Tilde{M}@>\pi>>\Tilde{M}$ is holomorphic.
\qed \edef\teoremace{\number\q.\number\x}\par\smallskip
%%%%%%%%%%%%%%%%%%%%%%%%%%%%%%%%%%%%%%%
We have\,:
\thm{Let $M$ be a generic $CR$ submanifold, of type $(n,k)$, of a complex
manifold $\Tilde{M}$. Let $x_0$ be a point of $M$ and 
$\xi_0\in \Dot{H}^0_{x_0}M$ a regular point of
$\Dot{H}^0M$. Then\,:
Then there exists an open neighborhood $\Tilde{\omega}$ of
$x_0$ in $\Tilde{M}$, an $(n-m)$-dimensional complex submanifold
$W$ of $\Tilde{\omega}$ with
$$x_0\in W\subset M\cap\Tilde{\omega}$$
and a  real valued smooth function $\rho:\Tilde{\omega}@>>>\Bbb{R}$
with $\rho(x)=0$ for $x\in M\cap\Tilde{\omega}$, 
$d^c\rho(x_0)=\xi_0$ and 
$\left.\left.
\dsize\frac{\partial\rho}{\partial z^\alpha}\right|_{\dsize W}\right.$
holomorphic in $W$ for $\alpha=1,\hdots,n+k$ (here  
$z^1,\hdots,z^{n+k}$ are holomorphic coordinates in
$\Tilde{\omega}$).}
\dimo 
It suffices to consider the situation where 
$\Tilde{M}=\Tilde{\omega}$ is a 
neighborhood of $x_0$ in $\Bbb{C}^{n+k}$. Then $T^*\Tilde{\omega}$ 
can be identified to the product manifold
$\Tilde{\omega}\times \left(\Bbb{C}^{n+k}\right)^*$, where
$\left(\Bbb{C}^{n+k}\right)^*$ is the space of $\Bbb{C}$-linear
forms in $\Bbb{C}^{n+k}$. If $M$ is described 
by (\defemme), then 
$H^0M$ gets identified to the real submanifold of 
$\Tilde{M}\times\left(\Bbb{C}^{n+k}\right)^*$\,:
$$\left.\left\{\left(z,(1/i)\sum_{j=1}^k{
a^j\partial\rho_j(z)}\right)\, 
\right|\, z\in M\, , \; a^1,\hdots,a^k\in\Bbb{R}\, \right\}\,.$$
By Lemma \teoremace, there is a complex submanifold
$V^*$ of dimension $(n-m)$ of $T^*\Tilde{M}$ that is contained
in $H^0M$, contains the point $\xi_0$,
and whose projection
$W=\pi(V^*)$ in $M$ is an $(n-m)$ smooth complex submanifold of
$\Tilde{M}$.
Define 
$\rho=\sum_{j=1}^k{a^j(z)\rho_j(z)}$ with real 
valued smooth
functions $a^j$ such that $\left(z,(1/i)\sum_{j=1}^k{
a^j(z)\partial\rho_j(z)}\right)$ belongs to $V^*$
when $z\in W$. \qed\edef\teoremab{\number\q.\number\x}
\medskip
%%%%%%%%%%%%%%%%%%%%%%%%%%%%%%%%%%%%%%%%%%%%%%%%%%%%%%%%%%%%%%%%%%%%%%%
%%%%%%%%%%%%%%%%%%%%%%%%%%%%%%%%%%%%%%%%%%%%%%%%%%%%%%%%%%%%%%%%%%%%%%%
\se{The main theorems (corrected version)}
\thm{Let $M$ be a locally embeddable $CR$ manifold of type $(n,k)$. 
Let $x_0\in M$ and assume that  
$\dot{H}^0_{x_0}M$ contains a regular point $\xi_0$ of $H^0M$ 
such that the Levi form $\Cal{L}_{\xi_0}$
has $q$ positive eigenvalues and $(n-q)$ 
eigenvalues which are $\leq 0$. Then the local cohomology groups
$\Cal H_{\bar\partial_M}^{q}\left((x_0),\Cal Q^{p,*}_M\right)=
\underset{U\ni x_0}\to{\varinjlim}
{H^q_{\bar\partial_M}\left(U,Q^{p,*}_M\right)}$
are infinite dimensional for all 
$0\leq p\leq n+k$.}\edef\teoda{\number\q.\number\x}
In fact a more general statement is valid. To formulate it, we consider
the $\bar\partial_M$-complex on currents (see [HN1], [NV]) and consider,
for an open subset $\omega$ of $M$, and $0\leq p\leq n+k$,
its
cohomology groups, that we denote by
 $H^q_{\bar\partial_M}\left(\omega,\Cal D'\otimes\Cal Q^{p,*}\right)$
(for $0\leq q\leq n$).
We define
the local cohomology on currents by:
 $$\Cal H^q_{\bar\partial_M}\left((x_0),\Cal D'\otimes\Cal Q^{p,*}\right)
=\underset_{\omega\ni x_0}\to{\varinjlim}{{
H^q_{\bar\partial_M}\left(\omega,\Cal D'\otimes Q^{p,*}_M\right)}}
\,
$$
 Note that we have a natural map 
 \form{\Cal H_{\bar\partial_M}^{q}\left((x_0),\Cal Q^{p,*}_M\right)
@>>>\Cal H^q_{\bar\partial_M}\left((x_0),\Cal D'\otimes\Cal Q^{p,*}\right)\,
.}
\edef\inclusio{\number\q.\number\t}
We have:
\thm{With the same assumptions of Theorem {\teoda},
the map (\inclusio) has an
infinite dimensional image.}\edef\teodb{\number\q.\number\x}
We give first the proof of Theorem {\teoda}, then indicate the small changes
needed to prove Theorem {\teodb}.\smallskip
\noindent
{\sc Proof of Theorem {\teoda}}\par
By [BHN1],  it suffices to show that 
$\Cal H_{\bar\partial_M}^{q}\left((x_0),\Cal Q^{p,*}_M\right)\neq 0$.
We can assume that $M$ is a generic $CR$ submanifold of
an open subset $\Tilde{M}$ of $\Bbb{C}^{n+k}$, and $M$ is defined
by (\defemme).
If $\Cal{L}_{\xi_0}$ has rank $n$ 
our statement reduces to \cite{AFN}. If  $\Cal{L}_{\xi_0}$ 
has rank $(n-d)$,
we utilize Theorem {\teoremab}\,.
By shrinking $\Tilde{M}$, and changing the holomorphic coordinates,
we put ourselves in the situation  where $x_0=0$, there is a 
$d$-dimensional complex linear space $W=\{z^{d+1}=0,\hdots,z^{n+k}=0\}$
such that
\form{x_0=0\in W\cap\Tilde{M}\subset M\subset\Tilde{M}}
and 
\form{\partial\rho_1(z)=idz^{n+1}\quad\roman{in}\quad W\cap\Tilde{M}\,,
\quad \partial\rho_j(0)=i\,{dz^{n+j}}|_0\quad\roman{for}\quad
j=1,\hdots,k\, .}\edef\formulacc{\number\q.\number\t}
 Set $\zeta=(z^1,\hdots,z^{n})$
and $t=(t^1,\hdots,t^k)$. 
By the implicit function theorem, after shrinking $\Tilde{M}$,
we can take the defining functions in the form\,:
\form{\rho_j(z)=s^j-h_j(z),\qquad\text{with}\quad
h_j(z)=h_j(\zeta,t)=0(2)\, ,}
and by ({\formulacc}) we obtain moreover that 
\form{h_1(z)=0\quad\roman{and}\quad
\dsize\frac{\partial h_1(z)}{\partial z^\alpha}
=0
\quad\roman{for}\quad z\in W\cap\Tilde{M}\, \quad\roman{and}\quad
\alpha=d+1,\hdots,n.}
By substituting the complex variable $z^{n+1}$ by 
\form{z^{n+1}-i
\sum_{\alpha,\beta=d+1}^{n+k}\frac{\partial^2h_1(0)}{\partial z^\alpha
\partial z^\beta} z^\alpha z^\beta\, ,}
and making a linear change of the variables $z^{d+1},\hdots,z^n$,
we obtain that\,:
\form{h_1(\zeta,t)=\sum_{\alpha=d+1}^{d+q}z^\alpha\bar z^\alpha\quad
-\;\sum_{\alpha=d+q+1}^n z^\alpha\bar z^\alpha\; +\; 0(3)\quad
\text{at $0$}\,.}\edef\formulacf{\number\q.\number\t}
By Taylor's formula we get\,: 
\form{\matrix\format\r&\l\\
h_1(\zeta,t)&=\dsize\sum_{\alpha,\beta=d+1}^{n+k}{\frac{\partial^2
h_1(z^1,\hdots,z^d,0,\hdots,0)}{\partial z^\alpha\partial\bar z^\beta}
\,z^\alpha\bar z^\beta}\,\\\vspace{2\jot}
& + \,\dsize
\Re\sum_{\alpha,\beta=d+1}^{n+k}{\frac{\partial^2
h_1(z^1,\hdots,z^d,0,\hdots,0)}{\partial z^\alpha\partial z^\beta}
\,z^\alpha z^\beta}
\\\vspace{2\jot}
\,&\quad +\,\dsize  o
\left(\sum_{\alpha,\beta=d+1}^n z^\alpha \bar z^\alpha\right)\, .
\endmatrix}
By (\formulacf) the second summand in the right hand side vanishes
at $0$. Thus we can find $r_0>0$ so that
$B(r_0)=\{|z|\leq r_0\}\Subset\Tilde{M}$ and 
\form{ \matrix\format\r&\l\\
h_1(z)
=h_1(\zeta,t)\quad &
\leq \quad 2\;\sum_{\alpha=d+1}^{d+q}z^\alpha\bar z^\alpha\quad
-\;(1/2)\sum_{\alpha=d+q+1}^n z^\alpha\bar z^\alpha\;\\\vspace{1\jot}
& \qquad\qquad\; +\sum_{j=1}^k \left(t^j\right)^2\\\vspace{2\jot}
&\qquad\qquad\;\qquad\qquad\roman{for}\quad
z\in B(r_0)\cap M\, .\endmatrix}
%%%%%%%%%%%%%%%%%%%%%%%%%%%%%%%%%%%%%%%%%%%%%%%%%%%%%%%%%%%%%%

Let us fix a  real number
 $\nu>2$ and set:
\form{\phi(z)=\phi(\zeta,t)=-it^1+
h_1(\zeta,t)-\nu\sum_{\alpha=d+1}^{d+q}{z^\alpha \bar z^\alpha}
-\nu\sum_{j=1}^k{(t^j+ih_j(\zeta,t))^2}\, .}
Since $h_j(z)=\Im z^{n+j}=0$ in $W\cap\Tilde{M}$, 
we can find a real number $r_1$ with $0<r_1\leq r_0$ such that
\form{\Re\phi(z)\leq 0\quad\roman{for}\quad z\in B(r_1)\cap M\, .}

%%%%%%%%%%%%%%
Define, for every real $\tau>0$
\form{f_\tau=e^{\left[
\frac{1}{\tau}\phi(z)\right]}\, 
dz^1\wedge\cdots\wedge dz^p \wedge d\bar z^{d+1}\wedge\cdots\wedge
 d\bar z^{d+q}\, .}
This is a smooth $(p,q)$-form, that defines a
form in $\Cal{Q}^{p,q}(B(r_1)\cap M)$ satisfying
\form{\bar\partial_M f_\tau=0\, \quad
\roman{in}\quad B(r_1)\cap M\,.}
We next define\,:
\form{\psi(z)=it^1-h_1(z)-\nu\sum_{\alpha=1}^d{z^\alpha 
\bar z^\alpha}-\nu\!\!\!\sum_{\alpha=d+q+1}^n{z^\alpha 
\bar z^\alpha}
-\nu\sum_{j=1}^k{(t^j+ih_j(\zeta,t))^2}\, .}
Then we can find a positive $r$ with $0<r\leq r_1\leq r_0$ such that
\form{\Re\psi(z)\leq -\frac{1}{2}\sum_{\alpha=1}^{n+k}{z^\alpha
\bar z^\alpha}\qquad\roman{for}\quad z\in B(r)
\cap M\, .}\edef\formulacq{\number\q.\number\t}
Now $\nu$ and $r$ are fixed and we set\,:
\form{g_\tau=e^{\left[
\frac{1}{\tau}\psi(z)\right]}\, 
dz^{p+1}\wedge\cdots\wedge dz^n\wedge\cdots\wedge dz^n\wedge d
\bar z^{1}\wedge\cdots d\bar z^d\wedge d\bar z^{d+q+1}
\wedge\cdots\wedge
 d\bar z^n\, .}
For each $\tau>0$ the form $g_\tau$ defines an element of
$\Cal{Q}^{n+k-p,n-q}(B(r)\cap M)$ and we have\,:
\form{
\bar\partial_M f_\tau=0\, ,\quad \bar\partial_M g_\tau=0\,
\quad\roman{in}\quad B(r)\cap M.}
Let $\chi=\chi(\zeta,t)$ denote a smooth real valued function defined in 
$\Bbb R^{2n+k}$
such that
$\chi=1$ for $|\zeta|^2+|t|^2<\frac{1}{2}$ and 
$\chi=0$ for $|\zeta|^2+|t|^2>\frac{2}{3}$.
If the Poincar\'e lemma is valid, we have,
for all $R>0$ sufficiently large, an a priori estimate:
\form{\matrix\format\l\\
\left|\dsize\int{\chi(R\zeta,Rt) f_\tau\wedge g_\tau}\right|\,
\\
{}
\\
\qquad\leq
\, C\left(
\underset{\smallmatrix
z\in B(r)\cap M\\
|a|+|b|\leq m
\endsmallmatrix}\to\sup
\left| D_z^a D_{\bar z}^b f_\tau\right|\right)\cdot\left(\sup\left|
R (\bar\partial_M\chi)(R\zeta,Rt) \wedge g_\tau\right|\right)
\endmatrix}\edef\aprioridiciotto{\number\q.\number\t}
with constants $C=C(R)$ and $m=m(R)$
which depend on $R$ but are independent of $\tau$.

\smallskip
Next we note that 
$$\phi(\zeta,t)+\psi(\zeta,t)
=-\nu\dsize\sum_{\alpha=1}^n{z^\alpha\bar{z}^\alpha}
-2\nu\dsize\sum_{j=1}^k{\left(t^j+ih_j(\zeta,t)\right)^2}\, .$$
Upon replacing $z$ by $z/\sqrt{\tau}$, we have:
$$\matrix\format\l\\
\dsize\int{\chi(R\zeta,Rt) f_\tau\wedge g_\tau}= (-1)^{q(n+k+d-p)}
\tau^{n+\frac{k}{2}}
\\
{}
\\
\qquad\qquad\times
\dsize\int{\chi\left(R\sqrt{\tau}(\zeta,t)\right)
\exp\left({{\smallmatrix
{-\nu\dsize\sum_{\alpha=1}^n{z^\alpha\bar{z}^\alpha}}\endsmallmatrix}
{\smallmatrix
{-2\nu\dsize\sum_{j=1}^k{\left(t^j+0(\sqrt{\tau})\right)^2}}
\endsmallmatrix}}\right)}\\
\qquad \qquad\qquad\qquad\qquad\qquad \times \;
dz^1\wedge\cdots\wedge dz^{n+k}
\wedge d\bar z^1\wedge\cdots\wedge d\bar z^n
\endmatrix$$
Therefore we obtain that
\form{
\CD
\tau^{-(n+\frac{k}{2})}
\left|\dsize\int{\chi(R\zeta,Rt) f_\tau\wedge g_\tau}\right|@>>> 
\roman{constant}>0
\endCD
\quad\text{for $\tau @>>>0$.}}\edef\formulaconlac{\number\q.\number\t}
\smallskip
Next we observe that
we have an estimate of the form:
\form{\underset{\smallmatrix
z\in B(r)\cap M\\
|a|+|b|\leq m
\endsmallmatrix}\to\sup
\left| D_z^a D_{\bar z}^b f_\tau\right|\leq c_1 \tau^{-m}}
%%%%%%%%%%%%%%%%%%%%%%
with a positive constant $c_1$ which is independent $\tau$,
because $\Re\phi\leq 0$ in $B(r)\cap M$.
\smallskip
On the other hand we have, for a positive constant $c_2$\,:
\form{\sup\left|
R (\bar\partial_M\chi)(Rz) \wedge g_\tau\right|\, \leq\,
c_2\cdot R\cdot \exp\left(-R^{-2}/(4\tau)\right)\, ,
}\edef\formulaventuno{\number\q.\number\t}
by (\formulacq).
\smallskip
By letting $\tau$ approach $0$, with any fixed large $R>0$,
we see that 
(\formulaconlac) cannot possibly hold true. 
This proves the theorem.
\medskip
\noindent
{\sc Proof of Theorem \teodb}\par
The only changes needed to obtain Theorem {\teodb} are in formulas 
(\aprioridiciotto)
and (\formulaventuno). 
In (\aprioridiciotto)
instead of the $\sup$ norm 
of $\bar\partial_Mg_\tau$ we need to introduce the $\sup$ norm of
its derivatives up to some finite order $m_1\geq 0$. 
This modifies (\formulaventuno) by a factor $\tau^{-m_1}$ in
the right hand side, but again the right hand side
of (\aprioridiciotto) tends to $0$ when $\tau@>>> 0$,
yielding a contradiction
that,
in view of
Theorem {\teoremafad}, proves the statement of Theorem \teodb.
\edef\sezionee{\number\q}
%%%%%%%%%%%%%%%%%%%%%%%%%%%%%%%%%%%%%%%%%%%%%%%%%%%%%%%%%%%%%%%%%%%%%%%%%%%
\se{Examples}
\noindent
 $1$.\quad Let $X$ be the complex manifold consisting of the
pairs $(L_1,L_3)$ where $L_i$ is a complex linear subspace of
 $\Bbb C^4$ of dimension $i$ and $L_1\subset L_3$. This is
a compact complex manifold of complex dimension $5$. 
Fix the canonical basis $e_1,e_2,e_3,e_4$ of $\Bbb C^4$. 
A local chart of $X$ in an open neighborhood $U$ of 
 $\left(\langle e_1\rangle,\langle e_1,e_2,e_3\rangle\right)$ in $X$ 
is given by the coordinates $z^1,\hdots,z^5$ where $L_1$ and $L_3$
are generated, respectively,  by the first column 
and by the three columns of the matrix:
$$\pmatrix
1&0&0\\
z^1&1&0\\
z^2&0&1\\
z^3&z^4&z^5
\endpmatrix\, .$$
Consider on $\Bbb C^2$ the structure of a right module over
the division ring $\Bbb H$ of the quaternions and 
let $M$ be the subset of $X$ consisting of
the pairs $(L_1,L_3)$ with $L_1\cdot\Bbb H\subset L_3$.
In $U$ the equation of $M$ is:
 $$ \left|\matrix
1&0&0&-\bar z^1\\
z^1&1&0&1\\
z^2&0&1&-\bar z^3\\
z^3&z^4&z^5&\bar z^2
\endmatrix\right|\, =\, 0$$
and therefore it is easy to verify that 
 $M$ is a $CR$ manifold of type $(3,2)$ and that  
the Levi form of $M$ has in every
non-zero characteristic codirection exactly one positive, one negative
and one zero eigenvalue. Thus at each point of $M$ the Poincar\'e
lemma fails in degree $1$.
\medskip
\noindent
 $2$.\quad Let $m_1,\hdots ,m_\ell$ be integers $\geq 1$ and fix complex
coordinates $w^1$, $w^2$, $z^1_j$, $\hdots$, $z^{m_j+1}_j$, 
$\zeta^1_j$, $\hdots$, $\zeta^{m_j}_j$, for $j=1,\hdots, \ell$
in $\Bbb C^{2(m_1+\cdots+m_\ell)+\ell+2}$.
We define a $CR$ submanifold of type 
 $(2[m_1+\cdots m_\ell]+\ell,2)$ of $\Bbb C^{2(m_1+\cdots+m_\ell)+\ell+2}$ by:
 $$M\, :=\,\cases
\Im w^1=\Im\dsize\sum_{j=1}^{\ell}
\dsize\sum_{h=1}^{m_j}{z_j^h\bar\zeta_j^h}\\
\Im w^2=\Im\dsize\sum_{j=1}^{\ell}\dsize\sum_{h=1}^{m_j}
{z_{j}^{h+1}\bar\zeta_j^h}\, .
\endcases$$
For every $\xi\neq 0$ in $H^0M$, the Levi form $\Cal{L}_{\xi}$
has $m_1+\cdots +m_\ell$ positive, $m_1+\cdots +m_\ell$ 
negative, and $\ell$ zero eigenvalues.
Thus the Poincar\'e lemma for $M$ is not valid at any point of $M$
in dimension $q=m_1+\cdots + m_\ell$. 
\medskip
\noindent
 $3$.\quad Let $N\subset\Bbb C^5$ be defined by
 $$N=\left\{(z,w)\in\Bbb C^2\times\Bbb C^3\, \left| \, 
{\matrix
\Im w^1=z^1\bar z^2+z^2\bar z^1\\
\Im w^2=z^1\bar z^3+z^3\bar z^1
\endmatrix}\right\}\right.$$
At each $\xi\in H^0M\setminus\{0\}$ the Levi form $\Cal{L}_{\xi}$
has one positive, one negative and one zero eigenvalue.
We set $t^j=\Re w^j$, $j=1,2$ and use $t^1, t^2, z^1, z^2, z^3$
as global coordinates on $N$.
The complex vector fields 
 $$\cases
\bar L_1=\frac{\partial}{\partial\bar z^1}-
iz^2\frac{\partial}{\partial\bar t^1}-
iz^3\frac{\partial}{\partial\bar t^2}\\
\bar L_2=\frac{\partial}{\partial\bar z^1}-
iz^1\frac{\partial}{\partial\bar t^1}\\
\bar L_3=\frac{\partial}{\partial\bar z^1}-
iz^1\frac{\partial}{\partial\bar t^2}
\endcases
 $$
give a basis of $T^{0,1}N$ at each point of $N$. 
According to our result, there exist smooth complex valued functions
 $\omega_1,\omega_2,\omega_3$, defined in a neighborhood $U$ of $0$
in $\Bbb R^2_t\times\Bbb C^3_z$ such that
 $\bar L_k\omega_j=\bar L_j\omega_k$ for all $1\leq j,k\leq 3$, but
the system $\bar L_ju=\omega_j$ for $j=1,2,3$ has no solution in a
neighborhood of $0$.
We define on $M=U\times C_{\zeta}$ a $CR$ structure by requiring that
the vector fields
 $$\cases
\bar L_1=\frac{\partial}{\partial\bar z^1}-
iz^2\frac{\partial}{\partial\bar t^1}-
iz^3\frac{\partial}{\partial\bar t^2}+\omega_1\frac{\partial}{\partial\zeta}
\\
\bar L_2=\frac{\partial}{\partial\bar z^1}-
iz^1\frac{\partial}{\partial\bar t^1}+\omega_2\frac{\partial}{\partial\zeta}
\\
\bar L_3=\frac{\partial}{\partial\bar z^1}-
iz^1\frac{\partial}{\partial\bar t^2}+\omega_3\frac{\partial}{\partial\zeta}\\
\bar L_4=\frac{\partial}{\partial\bar\zeta}
\endcases
 $$
form a basis of $T^{0,1}M$ at each point of $M$. This
 $M$ gives an example of a $CR$ manifold, of type $(4,2)$, such that the
Levi form has for each $\xi\in H^0M\setminus\{0\}$ 
at least $3$ eigenvalues $\geq 0$
and $3$ eigenvalues $\leq 0$ (actually one positive, one negative and
two zero eigenvalues), which according to [H1], [H2] is not locally
$CR$ embeddable near $0$. For similar examples see also [HN2]
\medskip
\noindent
 $4$.\quad 
Let $S^5$ be the unit sphere in $\Bbb{C}^3$ and consider the
manifold $M$ consisting of the complex lines that are tangent to
$S$. This is a $CR$ submanifold, of hypersurface type,
generically embedded into the complex $4$-dimensional manifold
$\Tilde{M}$ 
consisting of all complex lines of the 
three dimensional complex projective space $\Bbb{CP}^3$.
The $CR$ dimension of $M$ is $3$ and
its Levi form $\Cal{L}_\xi$ 
has at each nonzero $\xi\in H^0M$ 
one positive, one negative
and one zero eigenvalue. Thus 
for all $x_0\in M$ the Poincar\'e lemma is not valid at $x_0$
in degree $1$. 
\medskip
\noindent
$5$.\quad 
Consider the real hypersurface $M$ in $\Bbb{C}^3$ described by\,:
$$M=\{\Im z^3=z^1\bar z^1+(z^1+\bar z^1)^mz^2\bar z^2\,\}\,,$$
where $m$ is an integer $\geq 2$.
Then $\xi_0=\left(d\Re z^3\right)_{(0,0,0)}$ is a regular point.
Indeed, with $\rho=\Im z_3-z^1\bar z^1-(z^1+\bar z^1)^mz^2\bar z^2\,$,
we have $\partial\rho=(i/2)d z^3$ constant, and hence holomorphic,
along the holomorphic curve $W=\{z_1=0\,,\;z_3=0\}\subset M$.
Therefore $V=\{(z,\xi)\,|\, z\in W\,,\; \xi=d\Re z^3/2\}\subset
H^0M$ is a holomorphic curve in $T^*\Bbb{C}^3$. Since the Levi
form $\Cal{L}_{\xi_0}$ has rank $1$, we have in fact
$T_{\xi_0}V=H_{\xi_0}H^0M$ because
$T_{\xi_0}V\subset H_{\xi_0}H^0M$ and they have the same dimension.
The Levi form $\Cal{L}_{\xi_0}$ has one positive and one zero
eigenvalue. Then the Poincar\'e lemma fails in degree $1$ at
$x_0=(0,0,0)$.  Note that in this case the rank of the Levi form
is not maximal at the regular point $\xi_0$, and that the 
rank of the Levi form
$\Cal{L}_{\xi}$
is not constant  in any neighborhood of $\xi_0$ in $H^0M$.
\medskip
\noindent
$6$.\quad Let $\chi:\Bbb{R}@>>>\Bbb{R}$ be equal to $0$ for
$t\leq 1$ and equal to $t\exp(1/(1-t))$ when $t>1$. Then
$\chi$ is smooth and convex non decreasing. Let $q_1,q_2$ be positive
integers, with $q_1\geq 2$, $q_2\geq 1$,
and consider the hypersurface in $\Bbb{C}^{q_1+q_2}$\,:
$$M=\{(z,w)\in\Bbb{C}^{q_1}\times\Bbb{C}^{q_2}\, | \,
|z|^2+\chi(|w|^2)\,=\,1\,\}\, .$$
Then $M$ is the boundary of a smooth 
bounded convex set and therefore
all global cohomology groups of the tangential Cauchy-Riemann
complex
$\Cal{H}^j(M,\Cal{Q}^{p,*}_M)$ are zero for 
all $p=0,\hdots,q_1+q_2$ and $1\leq j\leq q_1+q_2-2$
(see e.g. \cite{N2}).
However, for all points $x_0=(z_0,w_0)\in M$ with $|w_0|<1$, if
$\xi_0\in H^0_{x_0}\setminus\{0\}$ the Levi form
$\Cal{L}_{\xi_0}$ has $q_1-1$ eigenvalues of the same sign and the
others equal to zero. Therefore the local cohomology groups
$\Cal{H}^{q_1-1}((x_0),\Cal{Q}^{p,*}_M)$ are infinite dimensional
for all $p=0,\hdots,q_1+q_2$.
\par

%%%%%%%%%%%%%%%%%%%%%%%%%%%%%%%%%%%%%%%%%%%%%%%%%%%%%%%%%%%%%%%%%%
%%%%%%%%%%%%%%%%%%%%%%%%%%%%%%%%%%%%%%%%%%%%%%%%%%%%%%%%%%%%%%%%%%
\se{Further remarks on the Poincar\'e Lemma}
We get back to the general situation considered in \S\sectfunctanal.
Assume that the complex (\complesso) admits the Poincar\'e lemma
in dimension $q$ 
at some point $x_0\in M$. Fix any Riemannian metric $\bold{g}$ on $M$ and
denote by $B(x_0,r)$ the ball of center $x_0$ and radius $r$ for the
distance defined by $\bold{g}$. By Theorem {\teorba}, for every $r>0$
there exists some $r'>0$ such that (\formulafai) is valid with
$\omega=B(x_0,r)$ and $\omega_1=B(x_0,r')$. For each $r>0$
denote by $\kappa_q(r)$ the supremum of these $r'>0$. In case 
the Poincar\'e lemma is not valid 
in dimension $q$ at $x_0$, and (\formulafai) does not hold for
$\omega=B(x_0,r)$ and any open $\emptyset\neq \omega_1\subset\omega$,
we set $\kappa_q(r)=0$.  We set\,:
\form{\nu_q^-(x_0)=\liminf_{r\searrow 0}\dsize\frac{\log{\kappa_q(r)}}{\log r}
\qquad\text{and}
\quad 
\nu_q^+(x_0)=\limsup_{r\searrow 0}\dsize\frac{\log{\kappa_q(r)}}{\log r}
\,.}
The values $\nu_q^-(x_0)\leq \nu_q^+(x_0)$ 
can be either real numbers $\geq 1$ or $+\infty$. [We are making
the convention that $\log 0=-\infty$.] If $\bold{d}$ is the distance
associated to the Riemannian metric $\bold{g}$ and
$\bold{d_1}$ the distance associated to another Riemannian metric
$\bold{g}_1$ on $M$, then there exist constants $\roman{C}_1>0$,
$\roman{C}_2>0$, $r_0>0$  such that\,:
\form{\bold{d}(x_0,x)\leq \roman{C}_1\bold{d}_1(x_0,x)
\leq \roman{C}_2\bold{d}(x_0,x) \quad\forall x\in M\quad
\roman{with}\quad \bold{d}(x_0,x)<r_0\, .}
Therefore we obtain\,:
\lem{The numbers $\nu_q^{\pm}(x_0)$ are 
independent of the Riemannian metric
$\bold{g}$\,.\par
If there exist a Riemannian metric $\bold g$ 
and a real $r_0>0$ such that\par\smallskip\centerline{
$H^q(B(x_0,r),\Cal{E}^*,A_*)=0$ for all $0<r<r_0$,}\par\smallskip
\noindent then $\nu_q^+(x_0)=\nu^-_q(x_0)=1$. }
The small balls of a Riemannian metric can be considered as
{\it convex sets}. Thus the condition that $\nu_q^+(x_0)>1$ has
the meaning that convexity is not sufficient for the vanishing
of the cohomology, but a small open subset (if there is any)
on which the cohomology vanishes in degree $q$ needs to have 
a special shape.
\par
  Suppose that $\nu_q^+(x_0)=1$.
  Then one could say that the cohomology
  {\it vanishes asymptotically} at $x_o$, in dimension $q$.
  This occurs for example if $\kappa_q(r)\simeq cr$, with $0<c<1$.
\par
       Suppose that the limit in $\nu_q^-(x_0)$ is infinite. Then
  one could say that the Poincar\'e lemma {\it fails asymptotically}
  at $x_o$, in dimension $q$.
         This occurs for example if 
$\kappa_q(r) \simeq  C\exp(-a/r)$, for
  positive constants $C$ and $a$.

\thm{Let $M$ be a locally embeddable $CR$ manifold of type $(n,k)$. 
Let $x_0\in M$ and assume that there exists $\xi\in H^0_{x_0}$ such that
the Levi form $\Cal{L}_{\xi}$ has $q$ positive eigenvalues 
and $(n-q)$ 
eigenvalues which are $\leq 0$. Let $\bold{g}$ be any Riemannian metric
on $M$. Then there are constants $r_0>0$ and $C>0$ such that, for
every $p=0,1,\hdots,n+k$, and $0<r'\leq r\leq r_0$, the maps\,:
\form{\Cal H_{\bar\partial_M}^{q}\left(B(x_0,r),\Cal Q^{p,*}_M\right)
@>>>
\Cal H_{\bar\partial_M}^{q}\left(B(x_0,r'),\Cal Q^{p,*}_M\right)
}\xdef\formulagc{\number\q.\number\t}
induced by the restriction have infinite dimensional image
if $r'>C r^{3/2}$. 
In particular $\nu_q^-(x_0)\geq 3/2$.}\edef\teorgc{\number\q.\number\x}

In fact a more general statement is valid, considering
the $\bar\partial_M$-complex on currents.
We have:
\thm{Let $M$ be a locally embeddable $CR$ manifold of type $(n,k)$. 
Let $x_0\in M$ and assume that there exists $\xi\in H^0_{x_0}$ such that
the Levi form $\Cal{L}_{\xi}$ has $q$ positive eigenvalues and $(n-q)$ 
eigenvalues which are $\leq 0$. Then 
there are constants $r_0>0$ and $C>0$ such that, for
every $p=0,1,\hdots,n+k$, and $0<r'\leq r\leq r_0$, the maps\,:
\form{\Cal H_{\bar\partial_M}^{q}\left(B(x_0,r),\Cal Q^{p,*}_M\right)
@>>>
\Cal H_{\bar\partial_M}^{q}\left(B(x_0,r'),
\Cal{D}'\otimes Q^{p,*}_M\right)}
 has an
infinite dimensional image
if $r'>C r^{3/2}$.}\edef\teorgd{\number\q.\number\x}
We give first the proof of Theorem \teorgc, then indicate the small changes
needed to prove Theorem \teorgd.\smallskip
\noindent
{\sc Proof of Theorem \teorgc}\par
Again by [BHN1], assuming as we can that $M$ is a generic $CR$ submanifold
of an open subset $\Tilde{M}$ of $\Bbb{C}^{n+k}$, we know that
  (\formulagc) has an infinite dimensional image
whenever it 
has a non zero image. Thus it will suffice to find conditions 
on $0<r'\leq r$
that are necessary in order that (\formulagc)
has zero image.\par
Therefore, we assume that $x_0=0$ and 
$M$ is a closed generic $CR$ submanifold of an 
open neighborhood $\Tilde{M}$ of $0$
in $\Bbb C^{n+k}$, where it is given by the equations:
\form{s^j=h_j(\zeta,t)=0(2)\qquad\text{for}\quad j=1,\hdots,k\, .}
Here 
$z^{n+j}=t^j+is^j$ with $t^j,s^j$ real for $j=1,\hdots,k$ and 
$\zeta=(z^1,\hdots,z^n)$.
\par
We can also assume that 
$$h_1(\zeta,t)=\sum_{\alpha=1}^n{\epsilon_\alpha z^\alpha \bar z^\alpha}+
\Im \sum_{\alpha=1}^n\sum_{j=1}^k{a_{\alpha,j}z^\alpha t^j}+
\sum_{j,\ell=1}^k{b_{j,\ell}t^jt^\ell}+0(3)\, ,$$
with $\epsilon_\alpha=1$ for $1\leq\alpha\leq q$ and
$\epsilon_\alpha\leq 0$ if $q<\alpha\leq n$.
By substituting the complex variable $z^{n+1}$ by
$z^{n+1}-\sum_{\alpha=1}^n\sum_{j=1}^k{a_{\alpha,j}z^\alpha z^{n+j}}
-i\sum_{j,\ell=1}^k{b_{j,\ell}z^{n+j}z^{n+\ell}}$
we can then reduce to the case where 
\form{h_1(\zeta,t)=\sum_{\alpha=1}^q{z^\alpha \bar z^\alpha}+
\sum_{\alpha=q+1}^n{\epsilon_\alpha z^\alpha \bar z^\alpha}
+0(3)\, .}
Let us fix a real number
 $\nu>1+\max_{q+1\leq \alpha\leq n}|\epsilon_\alpha|$ and set:
\form{\phi(z)=\phi(\zeta,t)
=-it^1+h_1(\zeta,t)-\nu\sum_{\alpha=1}^q{z^\alpha \bar z^\alpha}
-\nu\sum_{j=1}^k{(t^j+ih_j(\zeta,t))^2}}
\form{f_\tau=e^{\left[
\frac{1}{\tau}\phi(z)\right]}\, 
dz^1\wedge\cdots dz^p\wedge d\bar z^1\wedge\cdots\wedge
 d\bar z^q}
and 
\form{\psi(z)=
\psi(\zeta,t)=it^1-h_1(\zeta,t)-\nu\sum_{\alpha=q+1}^n{z^\alpha \bar z^\alpha}
-\nu\sum_{j=1}^k{(t^j+ih_j(\zeta,t))^2}}
\form{g_\tau=e^{\left[
\frac{1}{\tau}\psi(z)\right]}\, 
dz^{p+1}
\wedge\cdots\wedge dz^{n+k}\wedge d\bar z^{q+1}
\wedge\cdots\wedge
 d\bar z^n\, .}
Clearly we have
\form{
\bar\partial_M f_\tau=0\, ,\quad \bar\partial_M g_\tau=0\,.}
We use the same cutoff function $\chi$ as in \S{\sezionee}.
If (\formulagc) has zero image, 
for all $R>0$ with $R^{-1}<r'$, we have,  by Theorem {\teoremafab}, 
an a priori estimate:
\form{\matrix\format\l\\
\left|\dsize\int{\chi(R\zeta,Rt) f_\tau\wedge g_\tau}\right|\,
\\
{}
\\
\qquad\leq
\, C\left(
\underset{\smallmatrix
|\zeta|^2+|t|^2\leq r^2\\
|\alpha|\leq m
\endsmallmatrix}\to\sup
\left| D^\alpha f_\tau\right|\right)\cdot\left(\sup\left|
R (\bar\partial_M\chi)(R\zeta,Rt) \wedge g_\tau\right|\right)
\endmatrix}\edef\formgl{\number\q.\number\t}
with constants $C$ 
and $m$ which depends on $r,r'$ but are independent of $\tau$.

\smallskip
Next we note that 
$$\phi(\zeta,t)+\psi(\zeta,t)
=-\nu\dsize\sum_{\alpha=1}^n{z^\alpha\bar{z}^\alpha}
-2\nu\dsize\sum_{j=1}^k{\left(t^j+ih_j(\zeta,t)\right)^2}\, .$$
We have:
$$\matrix\format\l\\
\dsize\int{\chi(R\zeta,Rt) f_\tau\wedge g_\tau}\\\vspace{2\jot}
\; =
(-1)^{(n+k-p)q}
\tau^{n+\frac{k}{2}}\dsize\int{\chi\left(R\sqrt{\tau}(\zeta,t)\right)
\exp\left({{\smallmatrix
{-\nu\dsize\sum_{\alpha=1}^n{z^\alpha\bar{z}^\alpha}}\endsmallmatrix}
\!\!{\smallmatrix
{-2\nu\dsize\sum_{j=1}^k{\left(t^j+0(\sqrt{\tau})\right)^2}}
\endsmallmatrix}}\right)}\\
\qquad\qquad\qquad\qquad\qquad\qquad\qquad \times \;
dz^1\wedge\cdots\wedge dz^{n+k}
\wedge d\bar z^1\wedge\cdots\wedge d\bar z^n
\endmatrix$$
Therefore we obtain that
\form{
\CD
\tau^{-(n+\frac{k}{2})}
\left|\dsize\int{\chi(R\zeta,Rt) f_\tau\wedge g_\tau}\right|@>>> 
\roman{constant}>0
\endCD}
for $\tau @>>>0$.
\smallskip
Next we observe that
we have an estimate of the form:
\form{\underset{\smallmatrix
|\zeta|^2+|t|^2\leq r^2\\
|\alpha|\leq m
\endsmallmatrix}\to\sup
\left| D^\alpha f_\tau\right|\, \leq \,
\exp\left(c_1 r^3/\tau\right)}
with a positive constant $c_1$ which is independent of $\tau$.
Indeed the Taylor series of $\phi$ at $0$ has a purely imaginary first
degree term and a real second degree term which is $\leq 0$,
and the factor involving $\tau^{-m}$ is absorbed by the constant $c_1$.
\smallskip
On the other hand we have, for positive constants $c_2,c_3$:
\form{\sup\left|
R (\bar\partial_M\chi)(R\zeta,Rt) \wedge g_\tau\right|\, \leq\,
c_3\cdot R\cdot \exp\left(-c_2R^{-2}/\tau\right)\, .
}\edef\formgn{\number\q.\number\t}
Indeed the Taylor series of $\psi$ at $0$ has a purely imaginary
first order term and a real second order term which is negative
definite, while the form in the left hand side is different from zero
only for $|\zeta|^2+|t|^2>R^{-2}/2$.
\smallskip
Therefore 
$$c_2\, R^{-2}\leq c_1 r^3\quad\forall R>0\quad\roman{with}\quad
Rr'>1\,.$$
Hence $r'\leq \sqrt{c_1/c_2}\,\cdot\, r^{3/2}$.
This completes the proof of the theorem.
\medskip
\noindent
{\sc Proof of Theorem \teorgd}\par
The only changes needed to obtain Theorem {\teorgd} are in formulas 
(\formgl)
and (\formgn). Both remain valid when instead of the $\sup$ norm 
of $\bar\partial_Mg_\tau$ we need to introduce the $\sup$ norm of
its derivatives up to some finite order $m_1\geq 0$. In view of
Theorem  {\teoremafad}, 
we obtain the statement of Theorem \teorgd. 
\medskip\noindent
{\sc Example}\quad 
Let us consider, for an
integer $m\geq 2$, the $CR$ manifold $M$ of type $(2,1)$:
 $$M=\{(z^1,z^2,z^3)\in\Bbb C^3\, | \, |z^1|^2+
|z^2|^2+|z^3|^{2m}=1
\}
\, .$$
The Levi form is definite and nondegenerate at all points of $M$
where $z^3\neq 0$. But at the points of 
 $S=\{ |z^1|^2+|z^2|^2=1,\; z^3=0\}\subset M$
it is degenerate with one zero and one non zero eigenvalue.
The Poincar\'e lemma in degree $1$ is valid at all points
$x\in M\setminus S$ (see \cite{N2} or \cite{AH1},\cite{AH2}), 
but we have $\nu_1^-(x)\geq 3/2$
at all points $x$
of $S$. In particular the cohomology in degree $1$ of the intersection
of $M$ with a small Euclidean ball centered at $x\in S$ is
always infinite dimensional.
And the same is true if, instead
   of small Euclidean balls, one uses small balls in any
   Euclidean metric. Also, like in Example $6$, all of the global
   cohomology groups in degree $1$ are zero, because M is the
   smooth boundary of a weakly pseudoconvex domain (see \cite{BHN2}).

%%%%%%%%%%%%%%%%%%%%%%%%%%%%%%%%%%%%%%%%%%%%%%%%%%%%%%%%%%%%%%%%%%%%%%%%

\enddocument